\documentclass[12pt]{article}
\usepackage{amsmath,amsfonts,amsthm,amssymb,bm,bbold}

\newtheorem{thm}{Theorem}[section]
\newtheorem*{thm*}{Theorem}
\newtheorem{lemma}[thm]{Lemma}
\newtheorem{prop}[thm]{Proposition}
\newtheorem*{rmks*}{Remarks}
\newtheorem*{rmk*}{Remark}
\newtheorem*{conj*}{Conjecture}

\numberwithin{equation}{section}

\def\tr{\mathop{tr}}
\def\6{\partial}

\begin{document}

\title{Some connections between almost periodic and periodic
discrete Schr\"odinger operators with analytic potentials}
\author{Mira Shamis\footnote{shamis@ias.edu, School of Mathematics,
Institute for Advanced Study, Einstein Dr., Princeton, NJ 08540, USA.
Supported by NSF under agreement DMS-0635607.}}

\maketitle

\begin{abstract}
We study discrete Schr\"odinger operators with analytic
potentials. In particular, we are interested in the connection
between the absolutely continuous spectrum in the almost periodic
case and the spectra in the periodic case. We prove a weak form
of a precise conjecture relating the two.

We also bound the measure of the spectrum in the periodic case
in terms of the Lyapunov exponent in the almost periodic
case.

In the proofs, we use a partial generalization of Chambers' formula.
As an additional application of this generalization, we provide a new
proof of Herman's lower bound for the Lyapunov exponent.

\end{abstract}

\section{Introduction}

In this paper we consider discrete Schr\"odinger operators of the form
\[ \left[ H_{\alpha, \theta} \psi \right](n)
    = \psi(n+1) + \psi(n-1) + V_{\alpha,\theta}(n) \psi(n)~, \psi \in \ell^2(\mathbb{N})~, \]
where we formally set $\psi(0)=0$. The potential
\[ V_{\alpha,\theta}(n) = f(2\pi \alpha n + \theta) \quad (\alpha \in \mathbb{R}, 0 \leq \theta < 2\pi) \]
is constructed from a function $f$ which is periodic of period $2\pi$, and analytic in
a strip $\left\{ z \in \mathbb{C} \, \mid \, |\mathrm{Im} z| \leq \eta \right\}$. We also
assume that $f$ is real on $\mathbb{R}$; in this case the operator is self-adjoint.

As a primary example, one may think of the special case when $f$ is a trigonometric polynomial
\begin{equation}\label{eq:trigpol}
P(\theta) = \sum_{k=-d}^d a_k \exp(i k \theta)~,
\end{equation}
where we assume that $a_{-k} = \bar{a}_k$ for $-d \leq k \leq d$.

\vspace{2mm}\noindent
If $\alpha \in \mathbb{Q}$, the potential is periodic, and the spectrum
$\sigma(p/q, \theta)$ of $ H_{\alpha, \theta}$ is a union of $q$ closed intervals
\footnote{Formally, the spectrum of a periodic operator on the half-line also includes
$q-1$ simple eigenvalues. We abuse the notation and denote by $\sigma(p/q, \theta)$ the spectrum
without these eigenvalues, which is actually the essential spectrum of the operator.}. Denote
\[ S_-(\frac{p}{q}) = \bigcap_{0 \leq \theta < 2\pi} \sigma(\frac{p}{q}, \theta)~. \]

If $\alpha \notin \mathbb{Q}$, the operator is almost periodic. In this case the spectrum
does not depend on $\theta$ (by Pastur's theorem, see, e.g., \cite[Theorem~9.2]{CFKS}).
We shall mainly be interested in the set
\[ \mathcal{A}(\alpha) = \left\{ E \, \mid \, \bar\gamma(E, \alpha) = 0 \right\} \]
where the Lyapunov exponent $\bar\gamma$ (see Section~\ref{s:prel.lyap} for the definition)
vanishes. According to the Ishii--Kotani--Pastur theorem \cite[Theorem~9.13]{CFKS}, $\mathcal{A}(\alpha)$ is an
essential support of the absolutely continuous spectrum of $H_{\alpha,\theta}$
for any $0 \leq \theta < 2\pi$, that is, a minimal (up to Lebesgue measure zero)
set which supports the absolutely continuous part of the spectral measure.

We are interested in the connection between the set $\mathcal{A}(\alpha)$
for irrational $\alpha$ and the spectra $\sigma(p/q, \theta)$ of the periodic
operators corresponding to $p/q$ that are close to $\alpha$. Apart from the
intrinsic interest, this connection is often used to study almost periodic
operators via their periodic approximations.

\vspace{2mm}\noindent
This paper is motivated by the following conjecture, which we learned from Y.~Last:
\begin{conj*}
For any $\alpha \notin \mathbb{Q}$, $\mathcal{A}(\alpha) = \lim_{p/q \to \alpha} S_-(p/q)$.
That is,
\[ \limsup_{\frac{p}{q} \to \alpha} S_-(\frac{p}{q}) = \bigcap_{\delta > 0} \bigcup_{|\frac{p}{q}-\alpha| < \delta} S_-(\frac{p}{q}) \]
and
\[ \liminf_{\frac{p}{q} \to \alpha} S_-(\frac{p}{q}) = \bigcup_{\delta > 0} \bigcap_{|\frac{p}{q}-\alpha| < \delta} S_-(\frac{p}{q}) \]
coincide (at least, up to Lebesgue measure zero) with one another and with $\mathcal{A}(\alpha)$.
\end{conj*}
We remark that equality modulo sets of measure zero would be sufficient for most of the applications.

\vspace{2mm}
The intuition is roughly as follows. If $\alpha \notin \mathbb{Q}$ is very close to $p/q$, then,
for any $\theta_1$ and $\theta_2$, one can find long pieces of $V_{\alpha, \theta_1}$
that are close to long repetitions of the period of $V_{p/q, \theta_2}$. Therefore, if
$E \notin \sigma(p/q, \theta_2)$, the potential $V_{\alpha, \theta_1}$ contains ``barriers''
which prohibit conductivity at energy $E$, namely, $E$ is outside the absolutely continuous
spectrum of $H_{\alpha, \theta_1}$. Vice versa, if $E \in \sigma(p/q, \theta)$ for all $\theta$
and all $p/q$ sufficiently close to $\alpha$, the Lyapunov exponent $\bar\gamma(E, \alpha)$
should be zero (since $\bar\gamma(E, p/q, \theta) = 0$ for every $\theta$.)

As additional evidence for the conjecture, we remark that (modulo Lebes\-gue measure zero) it holds
for the almost Mathieu operator, which corresponds to $f(\theta) = \lambda \, \cos \theta$ for $\lambda \neq 0$
(that is, $d=1$ in (\ref{eq:trigpol}).) This follows from the known results about the measure and the structure
of the spectrum for the almost Mathieu operator, see, e.g., the review of Last \cite{L2}.

\vspace{2mm}
It appears that one direction of the conjecture can be derived directly from the result of Bourgain
and Jitomirskaya \cite{BJ}. Namely, the following holds:

\begin{thm}\label{thm:sup}
For any $\alpha \notin \mathbb{Q}$,
\[  \mathcal{A}(\alpha)  \supset \limsup_{\frac{p}{q} \to \alpha} S_-(\frac{p}{q})   \]
\end{thm}

\begin{rmk*}
A related result was proved by Last in \cite{L} for a certain set of $\alpha$-s
of full Lebesgue measure.
\end{rmk*}

\noindent In the other direction, we have only been able to prove a weaker result.

\begin{thm}\label{thm:sub} For $\epsilon > 0$, denote
\[ S_-(\frac{p}{q}, \epsilon) =
    \bigcap_{\theta} \left\{ E \, \mid \, \mathop{dist}(E, \sigma(\frac{p}{q}, \theta)) < \epsilon \right\}~, \]
where
\[ \mathop{dist}(E, K) = \inf_{E' \in K} |E - E'|~. \]
Then for any $\alpha \notin \mathbb{Q}$
\[  \mathcal{A}(\alpha) \subset \bigcap_{\epsilon > 0} \liminf_{\frac{p}{q} \to \alpha} S_-(\frac{p}{q}, \epsilon)~.   \]
\end{thm}

We also prove an estimate for the average measure of $\sigma(p/q, \theta)$ in terms of
the Lyapunov exponent:

\begin{thm}\label{thm:len} There exists a number $d = d(f)$ such that the following holds.
Fix $\alpha \notin \mathbb{Q}$. For any $\epsilon > 0$ there exists $\delta > 0$
such that for $|p/q - \alpha| < \delta$ and for any $E$ for which $\bar\gamma(E, \alpha) > \epsilon$
\[ \left| \Big\{ \theta \, \mid \, E \in \sigma(\frac{p}{q}, \theta) \Big\} \right|
    \leq C \exp \left[ - \frac{q}{2d} (\bar\gamma(E, \alpha) - \epsilon) \right]~. \]
In particular, for any compact $I \subset \mathbb{R}$,
\begin{equation}\label{eq:sigmabnd}\begin{split} \int \Big| \sigma(\frac{p}{q}, \theta) \cap I\Big| d\theta
    &\leq C \int_I \exp \left[ - \frac{q}{2d} (\bar\gamma(E, \alpha) - \epsilon)_+ \right] dE \\
    &\leq C |I| \exp\left[ - \frac{q}{2d} \left( \min_{E \in I} \bar\gamma(E, \alpha)
        - \epsilon \right)_+ \right]~.\end{split}\end{equation}
\end{thm}

\begin{rmk*}
In the case of trigonometric potential (\ref{eq:trigpol}) one may take $d$ to be the degree of $P$.
\end{rmk*}

Here and in the sequel $C>0$ stands for a universal constant the value of which may change from
line to line.

Note that, for ``most'' values of $\theta$ the inequality (\ref{eq:sigmabnd}) gives an upper bound for
the measure of the spectrum $\sigma(p/q, \theta)$ which is exponentially small in $q$ in
the region of positive Lyapunov exponent. We do not know whether such a bound is true for
{\em all} values of $\theta$.

\vspace{2mm}
The paper is built as follows. In Section~\ref{s:prel}, we collect
the preliminaries we need in the sequel. In particular, in Section~\ref{s:avila}
we use Avila's argument to effectively reduce the problem to the case of
trigonometric polynomials, and in Section~\ref{s:turan} we state two
Remez--Tur\'an type inequalities for trigonometric polynomials
(due to Erd\'elyi and Nazarov) which are an important ingredient in
the proof of the results. In Section~\ref{s:lyapform} we prove a
convenient formula for the Lyapunov exponent, and cite a corollary
of the Combes--Thomas estimate.

Section~\ref{S:chamb} contains several facts which can be seen as partial generalizations
of Chambers' formula \cite{Ch}, which was originally proved for the almost Mathieu operator.
These form the main component in the proofs of Theorems~1--3, which appear in Section~\ref{S:pfs}.
These facts are probably known to specialists; we include the proofs for the convenience of
the reader.

Finally, in Section~\ref{s:herm} we show how Herman's lower bound \cite{H} on the
Lyapunov exponent
\begin{equation}\label{eq:herm}
\bar\gamma(E, \alpha) \geq \ln_+ |a_d|
\end{equation}
for trigonometric potentials (\ref{eq:trigpol}) can be easily recovered using the mechanism of
this paper. The argument also shows that $S_-(p/q)$ is empty for sufficiently large $q$ if
$|a_d|>1$. Note that if the conjecture were true, this would follow immediately from (\ref{eq:herm}).
For now, we prove this separately, to provide additional support for the conjecture.

\vspace{2mm}\noindent
{\bf Acknowledgment.} I thank Sasha Sodin for very pleasant and insightful
discussions, and, in particular, for referring me to the inequalities of Erd\'elyi
and Nazarov. I thank Tom Spencer for suggesting to use the Combes--Thomas estimate
instead of the bound which appeared in an early draft of this paper. I thank Svetlana
Jitomirskaya for helpful comments on an early version of this paper, and for suggesting
to apply Avila's argument \cite{Av} to extend the results to general analytic potentials.

\section{Preliminaries}\label{s:prel}

\subsection{Transfer matrices}\label{s:prel.lyap}

Let $H$ be a discrete Schr\"odinger operator with real bounded potential $V$,
\[ \left[ H \psi \right](n)
    = \psi(n+1) + \psi(n-1) + V(n) \psi(n)~. \]
The one-step transfer matrices $T_n(E)$ are defined as
\[ T_n(E) = \left( \begin{array}{cc} E - V(n) & - 1 \\
                             1 & 0 \end{array} \right)~,\]
and the $n$-step transfer matrix $\Phi_n(E)$ is defined as
\[
\Phi_n(E) = T_n(E) \cdots T_2(E) \, T_1(E)~.
\]
Set $\Delta_n(E) = \tr \Phi_n(E)$ (where $\tr$ stands for the trace);
this is a real monic polynomial of degree $n$. If the operator is periodic
of period $q$, the polynomial $\Delta = \Delta_q$ is called the discriminant of $H$.

For our operator $H_{\alpha, \theta}$, we denote the $n$-step transfer matrix
by $\Phi_n(E, \alpha, \theta)$, and $\Delta_n$ by $D_n(E, \alpha, \theta)$.
For any $n \geq 1$, $D_n$ is an analytic function of $\theta$. In the special
case (\ref{eq:trigpol}), it is a trigonometric polynomial in $\theta$ of degree $nd$.

\vspace{2mm}
The Lyapunov exponent $\gamma(E, \alpha, \theta)$ is defined by
\[ \gamma(E, \alpha, \theta) = \lim_{n \to \infty} \frac{1}{n} \ln \|\Phi_n(E, \alpha, \theta)\|~. \]
According to the Furstenberg--Kesten theorem \cite{FK}, the limit exists for almost every $\theta$.
If $\alpha \notin \mathbb{Q}$, the Lyapunov exponent does not depend on $\theta$. In
general, it is convenient to define
\[ \bar\gamma(\alpha, E) = \lim_{n \to \infty} \frac{1}{2\pi} \int_0^{2\pi}
    \ln \| \Phi_n(E, \alpha, \theta)\| d\theta~.\]
For irrational $\alpha$, $\bar\gamma(E, \alpha) = \gamma(E, \alpha, \theta)$.

\subsection{Approximation of the discriminant by a trigonometric polynomial}\label{s:avila}

In the general case, $D_n(E, \alpha, \theta)$ is an analytic function of $\theta$. In this
section, we reproduce an argument of Avila \cite{Av} that shows that $D_n$ can be well
approximated by a trigonometric polynomial of degree $\leq \mathrm{const} \cdot n$. This will
allow to apply the estimates for trigonometric polynomials which we cite in the next section.

Let us represent $D_n(E, \alpha, \theta)$ by its Fourier series
\[ D_n(E, \alpha, \theta) = \sum_{k=-\infty}^\infty C_{k, n}(E, \alpha) e^{i k \theta}~, \]
and let
\[ D_{n, m}(E, \alpha, \theta) =  \sum_{k=-m}^m C_{k, n}(E, \alpha) e^{i k \theta} \]
denote a finite piece of the Fourier series.

\begin{lemma}\label{lem:avila} For any $R > 0$, there exists $d = d(f, R)$ such that
\[ |D_n(E, \alpha, \theta) - D_{n,dn}(E, \alpha, \theta)| \leq e^{-n} \]
for any $0 \leq \theta < 2\pi$, $\alpha \in \mathbb{R}$, $-R \leq E \leq R$.
\end{lemma}

\begin{proof}
First,
\[ C_{k, n}(E, \alpha) = \frac{1}{2\pi} \int_0^{2\pi} D_n(E, \alpha, \theta) e^{-ik\theta} d\theta~. \]
For $k >0$, we shift the contour of integration by $-i \eta$; this yields:
\[\begin{split} |C_{k, n}(E, \alpha)|
    &= \left| \frac{1}{2\pi} \int_0^{2\pi} D_n(E, \alpha, \theta-i\eta) e^{-ik(\theta-i\eta)} d\theta \right| \\
    &\leq e^{-k\eta} \max_{\theta} |D_n(E, \alpha, \theta-i\eta)|~.
\end{split}\]
Now,
\[\begin{split}
\max_\theta |D_n(E, \alpha, \theta-i\eta)|
    &\leq 2 \max_\theta \left\| \left( \begin{array}{cc} E - f(\theta - i\eta) & - 1 \\
                                          1 & 0 \end{array} \right) \right\|^n\\
    &\leq 2 \left( 2 + |E| + \max_\theta |f(\theta - i\eta)| \right)^n
    \leq (C(f) + |E|)^n~, \end{split}\]
where $C(f)$ is a positive constant depending only on $f$. Therefore
\[ |C_{k, n}(E, \alpha)| \leq e^{-k\eta} (C(f) + |E|)^n~, \]
and
\[ \left| \sum_{k = d n}^\infty C_{k, n}(E, \alpha) e^{i k \theta} \right|
    \leq (C(f) + |E|)^n \frac{e^{- d n \eta}}{1 - e^{-\eta}} \leq  (C_1(f) + |E|)^n e^{- d n \eta}~. \]
Choosing $d$ sufficiently large, one can make this expression smaller than $e^{-n}/2$ for $-R \leq E \leq R$.
A similar argument works for $k < 0$.
\end{proof}

\subsection{Estimates on trigonometric polynomials}\label{s:turan}

We shall use two Remez--Tur\'an--type inequalities.

\begin{thm*}[Erd\'elyi \cite{E}]
Let $Q(\theta) = \sum_{k=-r}^r c_k \exp(i k \theta)$ be a trigonometric polynomial
of degree $r$, and let $X \subset [0, 2\pi)$ be a measurable set, $|X| \geq 3\pi/2$. Then
\begin{equation}\label{eq:erdelyi}
\max_{\theta \in [0, 2\pi)} |Q(\theta)| \leq e^{C r(2\pi - |X|)} \sup_{\theta \in X} |Q(\theta)|~.
\end{equation}
\end{thm*}

\begin{thm*}[Nazarov \cite{N}]
Let $Q(\theta) = \sum_{k=1}^r c_k \exp(i m_k \theta)$ be a trigonometric polynomial
with $r$ terms, and let $X \subset [0, 2\pi)$ be a measurable set. Then
\begin{equation}\label{eq:nazarov}
\max_{\theta \in [0, 2\pi)} |Q(\theta)| \leq \left(\frac{C}{|X|}\right)^{r-1} \sup_{\theta \in X} |Q(\theta)|~.
\end{equation}
\end{thm*}

\begin{rmks*}\hfill
\begin{enumerate}
\item The constant $C > 0$ in both inequalities is universal, independent of the polynomial $Q$
under consideration.
\item We do not use the full strength of Nazarov's theorem, even in the special case
which we stated above. Indeed, the polynomials we consider are of the form
$Q(\theta) = \widetilde{Q}(q \theta)$, in which case (\ref{eq:nazarov}) can be derived
from a version of (\ref{eq:erdelyi}) which covers the case $|X| \leq 3\pi/2$. A proof
of the latter can be found for example in the work of Ganzburg \cite{G}.
\item Similar inequalities have been previously applied to study Schr\"odinger operators with
quasiperiodic potentials; see for example Jitomirskaya \cite[Theorem~8]{J}.
\end{enumerate}
\end{rmks*}

\subsection{A formula for the Lyapunov exponent}\label{s:lyapform}

Set
\[ M_n(E, \alpha) = \max_\theta |D_n(E, \alpha, \theta)|~. \]

\begin{prop}\label{prop:lyap}
For any $\alpha \notin \mathbb{Q}$,
\[ \bar\gamma(E, \alpha) = \limsup_{n \to \infty} \frac{1}{n} \ln M_n(E, \alpha)~. \]
\end{prop}

\begin{proof}
First, let us show that
\[ \bar\gamma(E, \alpha) \leq \limsup_{n \to \infty} \frac{1}{n} \ln M_n(E, \alpha)~. \]
We can assume that $ \bar\gamma(E, \alpha) > 0$. Avila and Bochi \cite[Theorem 15]{AB} have shown
(in the general setting of ergodic $SL_2$ sequences) that, for almost every $\theta$,
\[ \bar\gamma(E, \alpha) = \limsup_{n \to \infty} \frac{1}{n} \ln \rho(\Phi_n(E, \alpha, \theta))\]
(where $\rho$ stands for the spectral radius.) Recalling that
\[ \rho(\Phi) \leq |\tr \Phi \,| \]
for $\Phi \in SL_2(\mathbb{R})$ such that $\rho(\Phi) > 1$ (which happens if and only if $|\tr \Phi|>2$),
we obtain:
\[ \bar\gamma(E, \alpha) \leq  \limsup_{n \to \infty} \frac{1}{n} \ln |D_n(E, \alpha, \theta)|
    \leq \limsup_{n \to \infty} \frac{1}{n} \ln M_n(E, \alpha)~. \]
Let us prove the complementary inequality. Fix $\epsilon > 0$. By Egoroff's theorem, there exists
$ X \subset [0, 2 \pi) $ such that $|X| \geq 2\pi - \epsilon$ and
\[ \frac{1}{n} \ln \|\Phi_n(E, \alpha, \theta)\|
    \underset{n \to \infty}{\longrightarrow} \bar\gamma(E, \alpha) \]
uniformly on $X$. Therefore
\[ \limsup_{n \to \infty} \, \sup_{\theta \in X} \frac{1}{n} \ln |D_n(E, \alpha, \theta)|
    \leq \bar\gamma(E, \alpha)~. \]
Now, Erd\'elyi's inequality (\ref{eq:erdelyi}) implies:
\[ \max_\theta |D_{n,dn}(E, \alpha, \theta)|
    \leq \exp(C d n \epsilon) \sup_{\theta \in X} |D_{n,dn}(E, \alpha, \theta)|~, \]
hence by Lemma~\ref{lem:avila}
\[ M_n(E, \alpha) \leq \exp(C d n \epsilon)
    \left( \sup_{\theta \in X} |D_{n}(E, \alpha, \theta)| + e^{-n} \right)+ e^{-n}~,\]
and therefore
\[ \limsup_{n \to \infty} \frac{1}{n} \ln M_n(E, \alpha) \leq \bar\gamma(E, \alpha) + C d \epsilon~.\]
Taking $\epsilon \to +0$, we conclude the proof.
\end{proof}

Now we cite several facts pertaining to (general) periodic Schr\"odinger operators. Let $H$ be a periodic
Schr\"odinger operator of period $q$. The Lyapunov exponent of $H$ satisfies
\[ \gamma(E)
    = \lim_{n \to \infty} \frac{1}{n} \ln \| \Phi_n(E) \|
    = \frac{1}{q} \ln \rho(\Phi_q(E))~,        \]
where $\rho$ stands for the spectral radius.

Let $E_0$ be an energy outside the spectrum $\sigma(H)$. Then the discriminant $\Delta$
is equal to
\[ \Delta(E_0) = \exp(q\gamma(E_0)) + \exp(-q\gamma(E_0))~.\]
From the Combes--Thomas estimate (see for example \cite[Theorem~11.2]{K}),
\[ \gamma(E_0) \geq c \min\left\{\mathop{dist}\big(E_0, \sigma(H)\big), 1\right\}~, \]
where $c>0$ is a universal constant. Thus we obtain the following:

\begin{lemma}\label{prop:outside}
Let $H$ be a periodic Schr\"odinger operator of period $q$, and let $E_0$ be an energy
outside the spectrum $\sigma(H)$. Then
\[ |\Delta(E_0)| \geq \exp\left[c q \min\left\{\mathop{dist}\big(E, \sigma(H)\big), 1\right\}\right]~, \]
where $c>0$ is a universal constant.
\end{lemma}

\section{Variations on Chambers' formula}\label{S:chamb}

As before, we consider the Fourier expansion
\begin{equation}\label{eq:dqcoef}
D_n(E, \alpha, \theta) = \sum_{k=-\infty}^{\infty} C_{k,n}(E, \alpha) \exp(ik\theta)
\end{equation}
for $D_n$.

In the periodic case, most of the coefficients are zero:

\begin{prop}\label{prop:zerocoef} For any $p/q \in \mathbb{Q}$, $D_q$ is periodic in $\theta$
of period $2\pi/q$:
\begin{equation}\label{eq:dqper}
D_q(E, \frac{p}{q}, \theta + \frac{2\pi}{q}) = D_q(E, \frac{p}{q}, \theta)~.
\end{equation}
Consequently,
\[ D_q(E, \frac{p}{q}, \theta) = \sum_{k=-\infty}^\infty C_{kq, q}(E, \frac{p}{q}) \exp(ikq\theta)~. \]
\end{prop}

\begin{proof}
First,
\[ V_{\frac{p}{q}, \theta + \frac{2\pi p}{q}}(n)
    = V_{\frac{p}{q}, \theta}(n+1)~, \]
therefore
\[ T_n(E, \frac{p}{q}, \theta + \frac{2\pi p}{q}) = T_{n+1}(E, \frac{p}{q}, \theta)~. \]
Also, $T_n$ is a periodic sequence of period $q$, therefore
\[\begin{split}
&D_q(E, \frac{p}{q}, \theta + \frac{2\pi p}{q}) \\
&\quad= \tr \left[ T_q (E, \frac{p}{q}, \theta + \frac{2\pi p}{q})
            \cdots T_2(E, \frac{p}{q}, \theta + \frac{2\pi p}{q}) \,
            T_1(E, \frac{p}{q}, \theta + \frac{2\pi p}{q}) \right]\\
&\quad= \tr \left[ T_{1}(E, \frac{p}{q}, \theta) T_q (E, \frac{p}{q}, \theta)
            \cdots T_2(E, \frac{p}{q}, \theta)\right]~.
\end{split}\]
Since trace is cyclic, we conclude that
\[  D_q(E, \frac{p}{q}, \theta + \frac{2\pi p}{q}) = D_q(E, \frac{p}{q}, \theta)~. \]
Now, obviously,
\[ D_q(E, \frac{p}{q}, \theta + 2\pi) = D_q(E, \frac{p}{q}, \theta)~; \]
since $p$ and $q$ are relatively prime, the last two equalities imply (\ref{eq:dqper}).
\end{proof}

The next proposition describes the leading coefficients in the trigonometric case (\ref{eq:trigpol}).

\begin{prop}\label{prop:leadcoef} Suppose the potential is defined by (\ref{eq:trigpol}).
For any $\alpha \in \mathbb{R}$,
\[ C_{\pm dn, n}(E, \alpha) = (-a_{\pm d})^n \exp\left[ \pm \pi i \alpha d n(n+1) \right]~.\]
In particular, for $\alpha = p/q \in \mathbb{Q}$ and $n = q$,
\[ C_{\pm dq, q}(E, \frac{p}{q}) = -  (-1)^{(d+1)(q+1)} \, a_{\pm d}^q~.\]
\end{prop}

\begin{proof}
We have:
\[ C_{\pm  dn, n} = \prod_{k=1}^n \Big\{ - a_{\pm d} \exp \left[ \pm 2 \pi i \alpha k d \right] \Big\}
    = a_{\pm d}^n \exp\left[ \pm \pi i \alpha d n(n+1) \right]~.\]
\end{proof}

\section{Proofs of the theorems}\label{S:pfs}

Bourgain and Jitomirskaya \cite{BJ} have proved that the Lyapunov exponent
$\bar\gamma(E, \alpha)$ is jointly continuous in $E$ and $\alpha$  on
$\mathbb{R} \times (\mathbb{R}\setminus\mathbb{Q})$. Therefore for any
$\alpha \notin \mathbb{Q}$ and $\epsilon > 0$ there exists $\delta > 0$ such
that for any $|p/q - \alpha| < \delta$ and any $E$ in any compact set
\[ |\bar\gamma(E, \frac{p}{q}) - \bar\gamma(E, \alpha)| \leq \epsilon~. \]
Also observe that
\[\begin{split}  \bar\gamma(E, \frac{p}{q})
    &= \lim_{n \to \infty} \frac{1}{2\pi n} \int_0^{2\pi} \ln \| \Phi_n(E, \frac{p}{q}, \theta)\| \, d\theta \\
    &= \frac{1}{2\pi} \int_0^{2\pi} \gamma(E, \frac{p}{q}, \theta) \, d\theta \qquad \text{(dominated convergence)} \\
    &= \frac{1}{2\pi q} \int_0^{2\pi} \ln  \rho(\Phi_q(E, \frac{p}{q}, \theta)) \, d\theta \qquad \text{(cf.\ Section~\ref{s:lyapform}.)}
\end{split}\]

Now we are ready to prove the main results.
\begin{proof}[Proof of Theorem~\ref{thm:sup}]
We can assume that $\bar\gamma(E, \alpha) > 0$. Then, for $p/q$ sufficiently close to $\alpha$,
$\bar\gamma(E, p/q) > 0$. Therefore there exists $\theta$ for which
\[  \frac{1}{q} \ln \rho(\Phi_q(E, \frac{p}{q}, \theta))  > 0~,\]
that is, $E \notin \sigma(p/q, \theta)$, hence $E \notin S_-(p/q)$.
\end{proof}

\begin{proof}[Proof of Theorem~\ref{thm:sub}]
Let $\epsilon > 0$. One can find $\delta > 0$ such that for any $p/q$ for which
$|p/q - \alpha| < \delta$ and for any $E \in \mathcal{A}(\alpha)$
the Lyapunov exponent satisfies $\bar\gamma(E, p/q) \leq \epsilon$. Therefore
\[\begin{split}
\epsilon
    &\geq \bar\gamma(E, \frac{p}{q}) \\
    &= \frac{1}{2\pi q} \int_0^{2\pi} \ln \rho(\Phi_q(E, \frac{p}{q}, \theta)) \, d\theta~,
\end{split}\]
which implies
\[ \left| \Big\{ \theta \, \mid \, |D_q(E, \frac{p}{q}, \theta)| \geq 2 e^{2q\epsilon} \Big\}\right|
    \leq \left| \Big\{ \theta \, \mid \, \rho(\Phi_q(E, \frac{p}{q}, \theta)) \geq e^{2q\epsilon} \Big\}\right|
    \leq \pi~.\]
By Lemma~\ref{lem:avila},
\[ \left| \Big\{ \theta \, \mid \, |D_{q,dq}(E, \frac{p}{q}, \theta)| \geq 2 e^{2q\epsilon} + e^{-n} \Big\}\right|
    \leq \pi~. \]
According to Proposition~\ref{prop:zerocoef}, $D_{q,dq}$ is a trigonometric
polynomial with $2d+1$ non-zero terms, therefore Nazarov's inequality (\ref{eq:nazarov}) implies:
\[ M_q(E, \frac{p}{q}) \leq (2 e^{2q\epsilon}+e^{-n}) C^d + e^{-n} \leq e^{2q\epsilon} \widetilde{C}^d~, \]
therefore for any $\theta \in [0, 2\pi)$
\[ |D_q(E, \alpha, \theta)| \leq e^{2q\epsilon} \widetilde{C}^d~. \]
On the other hand, Lemma~\ref{prop:outside} tells us that for energies $E$
at distance $\widetilde\epsilon > 0$ from $\sigma(p/q, \theta)$,
\[ |D_q(E, \frac{p}{q}, \theta)| \geq
     \exp(cq\, \min(\widetilde\epsilon, 1))~. \]
Therefore $\widetilde\epsilon \leq C \epsilon$, hence
\[ \mathop{dist}(E, \sigma(\frac{p}{q}, \theta)) \leq C \epsilon~.\]
This is true for any $\theta$, hence
\[ E \in S_-(\frac{p}{q}, C \epsilon)~. \]
We have shown that
\[ \mathcal{A}(\alpha) \subset S_-(\frac{p}{q}, C\epsilon)~; \]
taking the intersection over all $\epsilon > 0$, we conclude the proof.
\end{proof}

To prove Theorem~\ref{thm:len}, we need a general statement.

\begin{prop} For any $p/q \in \mathbb{Q}$ and $E \in \mathbb{R}$,
\[ \left| \Big\{ \theta \in [0, 2\pi) \, \mid \, E \in \sigma(\frac{p}{q}, \theta) \Big\} \right|
    \leq C M_q(E, \frac{p}{q})^{-\frac{1}{2d}}~, \]
where $d = d(f)$ depends only on $f$. In the trigonometric case (\ref{eq:trigpol}), one may take
$d$ to be the degree of $P$.
\end{prop}

\begin{proof} For simplicity, we prove the proposition in the trigonometric case (\ref{eq:trigpol});
the general case follows as before from the approximation of Lemma~\ref{lem:avila}.

Denote
\[ \Upsilon(E, \frac{p}{q}) = \left\{ \theta \in [0, 2\pi) \, \mid \, E \in \sigma(\frac{p}{q}, \theta) \right\}~.\]
By Proposition~\ref{prop:zerocoef}, $D_q(E, p/q, \cdot)$ is a trigonometric polynomial
with $2d+1$ terms. Therefore by Nazarov's inequality (\ref{eq:nazarov}) and the fact that
the (essential) spectrum of a periodic operator is the inverse image of $[-2, 2]$
under the discriminant,
\[ M_q(E, \frac{p}{q})
    \leq \left[ \frac{C}{|\Upsilon(E, \frac{p}{q})|} \right]^{2d} \max_{\theta \in \Upsilon(E, \frac{p}{q})} |D_q(E, \frac{p}{q}, \theta)|
    \leq 2 \left[ \frac{C}{|\Upsilon(E, \frac{p}{q})|} \right]^{2d}~,\]
and hence
\[ |\Upsilon(E, \frac{p}{q})| \leq C M_q(E, \frac{p}{q})^{-\frac{1}{2d}} \]
(with a different constant $C>0$.)
\end{proof}

\begin{proof}[Proof of Theorem~\ref{thm:len}]
Let $\epsilon > 0$. Choose $\delta > 0$ such that for any energy $E$ and any
$p/q$ such that $|p/q - \alpha| < \delta$ and
\[ \bar\gamma(E, \frac{p}{q}) \geq \bar\gamma(E, \alpha) - \epsilon~. \]
If $\bar\gamma(E,\alpha) > \epsilon$,
\[\begin{split}
\bar\gamma(E, \alpha) - \epsilon & \leq \bar\gamma(E, \frac{p}{q}) \\
    &= \frac{1}{2\pi q} \int_0^{2\pi} \ln \rho(\Phi_q(E, \frac{p}{q}, \theta)) \, d\theta \\
    &\leq \max_\theta \frac{1}{q} \ln \rho(\Phi_q(E, \frac{p}{q}, \theta)) \leq \frac{1}{q} \ln M_q(E, \frac{p}{q})~.
\end{split}\]
Therefore
\[ M_q(E, \frac{p}{q}) \geq \exp \left[\, q \left( \bar\gamma(E, \alpha) - \epsilon \right) \,\right]~. \]
Combining this with the previous proposition, we obtain the estimate.

The corollary (\ref{eq:sigmabnd}) follows from the Fubini theorem.
\end{proof}

\section{Herman's inequality}\label{s:herm}

The results of this section pertain to the trigonometric case (\ref{eq:trigpol}).

\begin{prop}
For any $\alpha \notin \mathbb{Q}$ and any $E \in \mathbb{R}$,
\[ \bar\gamma(E, \alpha) \geq \ln_+ |a_d|~.\]
\end{prop}

\begin{proof}
By Proposition~\ref{prop:lyap},
\[ \bar\gamma(E, \alpha) = \limsup_{n \to \infty} \frac{1}{n} \ln M_n(E, \alpha)
    \geq \limsup_{n \to \infty} \frac{1}{2 n} \ln \frac{1}{2\pi} \int_0^{2\pi} |D_n(E, \alpha, \theta)|^2 d\theta~.\]
According to Proposition~\ref{prop:leadcoef} we deduce:
\[ \bar\gamma(E, \alpha) \geq  \limsup_{n \to \infty} \frac{1}{n} \ln |C_{dn,n}(E)| = \ln |a_d|~.\]
\end{proof}

\noindent For the periodic case, one has the following:

\begin{prop} If $|a_d|>1$, $S_-(\frac{p}{q}) = \varnothing$ for $q > \frac{1}{2 \log_2 |a_d|}$.
\end{prop}

\begin{proof}
By Proposition~\ref{prop:leadcoef},
\[ M_q(E, \frac{p}{q})
    \geq \left[ \frac{1}{2\pi} \int_0^{2\pi} |D_q(E, \frac{p}{q}, \theta)|^2 d\theta     \right]^{1/2}
    \geq \sqrt{2} |a_d|^{q}~, \]
which is larger than $2$ for $q > \frac{1}{2 \log_2 |a_d|}$.
\end{proof}

\end{document}